\documentclass[11pt]{amsart} 
\usepackage[cp1251]{inputenc}
\usepackage[english]{babel}
\usepackage{amsmath}
\usepackage{amssymb}
\usepackage{amsfonts}
\usepackage{graphicx}[]

\newtheorem{theorem}{Theorem}[section]
\newtheorem{definition}{Definition}[section]
\newtheorem{corollary}{Corollary}[section]
\newtheorem{proposition}{Proposition}[section]
\newtheorem{comment}{Comment}[section]

	\title 
	[Decompositions of Finitely Additive Markov Chains]
	{Decompositions of finitely additive \\ Markov chains and invariant measures\\
		in discrete space}
	\author{{A.I. ZHDANOK}, {A.K. KHURUMA}}%
	\address{Alexander Ivanovich Zhdanok  
\newline\hphantom{A.I. ZHDANOK} Institute for Information Transmission Problems RAS, 
\newline\hphantom{A.I. ZHDANOK} Moscow, Russia}%
\email{zhdanok@inbox.ru}%

\address{Anna Kys-oolovna Khuruma  
	\newline\hphantom{A.K. KHURUMA} Tuvan State University, Kyzyl, Russia}
\email{huruma@list.ru}%

	\thanks{\rm 
	{\bf Acknowledgments:} 
	The reported study was funded by the Russian Foundation for Basic Research (both authors), project number 20-01-00575-a.
	}
		
\begin{document} 

\maketitle {\small
\begin{quote}
\noindent{\sc Abstract. }
 In this paper, we consider general Markov chains (MC), specified by the transition probability (kernel) $ P (x, E) $, finitely additive in the second argument.
 Such MC are studied within the framework of the functional operator treatment.
  The state space (phase space) of the MC $ X $ has any cardinality, and the sigma-algebra $ \Sigma $ is discrete, i.e. is the set of all subsets in $ X $.
   This construction of the phase space $ (X, \Sigma) $ allows us to decompose the Markov kernel $ P (x, E) $ into the sum of two components - countably additive and purely finitely additive in the second argument and measurable in the first argument.
   It is shown that the countably additive kernel  is atomic.
Some properties of Markov operators with a purely finitely additive kernel and their invariant measures are studied.
 A class of combined finitely additive MC and two of its subclasses are introduced, and some properties of their invariant measures are proved. Some asymptotic regularities of such MC are revealed.

{\bf Key words:} finitely additive measures, Markov chains, Markov operators, invariant measures.

MSC: 60J05, 46E27.
 \end{quote}
}

\section {\bf Introduction}


       Classical Markov chains (MC) are interpreted as random Markov processes
        with discrete time (in the usual sense) in the  phase space $(X, \Sigma)$, where $ X $ is some set (space), and $\Sigma$ is some sigma-algebra of subsets in $X$. We consider time-homogeneous MCs.    
 
      If $X$ is an arbitrary infinite set, without separating any structure in it, except for the sigma-algebra $\Sigma$, then such MCs are called general.
       
        There are two approaches to studying MC. Historically, the first approach is to describe the MC in the language of random variables. However, if the MC is general, and there are no linear operations in $ X $, then this approach is not applicable. It is also difficult to apply when finitely additive measures are used in the study of MCs.

     In 1937 N.M. Krylov and N.N. Bogolyubov proposed an operator approach to the study of general MC, which was then developed in detail in the article by K. Yosida and S. Kakutani \cite{YoKa:1}. Its essence is that the MC is given by the transition function (probability) $P(x, E), x \in X, E \in \Sigma$, which, as a kernel, defines two dual integral Markov operators $ T $ and
        $ A $ in spaces of measurable functions and in spaces of measures, respectively.  
   
      A Markov chain is identified with an iterative sequence of probability measures generated by the second Markov operator $\mu_{n}=A\mu_{n-1}=A^{n} \mu_{0}$ with an arbitrary initial measure $ \mu_{0} $.
      We use this approach in this work. \\
     
       In the classical operator theory of MC, the transition function (probability) $ P (x,\cdot) $ is assumed to be a countably additive measure with respect to the second argument. At the same time, in the economic theory of games, developed since the 1960s by L. Dubbins L. Savage \cite{DubSav1} and their numerous students and followers, it became necessary to involve in the construction and study of specific random processes and finitely additive probability measures. In particular, in the paper \cite {PS1}, some constructions of finitely additive measures close to Markov chains were constructed and investigated.

     In 1981, S. Ramakrishnan \cite{Rama1}, based on the work of \cite {PS1}, developed in the language of strategies a corresponding construction of a new object, which he called finitely additive Markov chains. These chains are spawned  a transition function (strategy), finitely additive in the second argument. The phase space $ (X, \Sigma) $ in \cite {Rama1} is a discrete set with the sigma-algebra of all its subsets. Within the framework of this construction, a number of theorems on the properties of these specific chains in terms of game theory are proved in \cite {Rama1}. Some additional questions on this topic are also discussed in further publications of Ramakrishnan (see, for example, \cite {Rama2}).
       
       Note that in the same 1981, the article \cite {Zhd00} of one of the co-authors of this paper A. Zhdanok was published simultaneously and independently, in which Markov sequences of finitely additive measures were studied, but generated by a Markov operator with a countably additive transition function (kernel). Within the framework of this approach, a number of ergodic-type theorems were further obtained for MC in general and metric phase space using finitely additive measures (see \cite {Zhd01} and \cite {Zhd02}). The works of other authors in this direction are also indicated there.
            
       Markov chains generated by a finitely additive transition function are also considered in our papers  A. Zhdanok \cite {Zhd01}, and A. Huruma \cite {XAK1}, \cite {XAK2}. Corresponding comments will be given in the text of the article.\\
       
       In this paper, we study Markov chains generated by a transition function, finitely additive in the second argument, given on a discrete space. However, it does not use any specific features of the theory of games, and a different range of problems is solved.
       
         We construct and study the decomposition of finitely additive Markov chains (i.e. Markov kernels) on a discrete set $ X $ of any cardinality (from an infinite power scale and when accepting the continuum hypothesis) into countably additive and purely finitely additive components and clarify the question of whether they have invariant measures of various types. We introduce and study the class of combined Markov chains. \\
    
        However, at the beginning, in the next item 2, we give an operator construction of a general Markov chain with a countably additive transition function on an arbitrary measurable space, specifying some details we need. We will use this construction further for finitely additive transition functions. The main new results obtained are presented and proved in items 3 and 4.
       
\section {\bf Definitions, notation and some information} 

Let $X$ be an arbitrary infinite set and $\Sigma$ the sigma-algebra of its subsets containing all one-point subsets from $X$. Let $B(X,\Sigma)$ denote the Banach space of bounded $\Sigma$ -measurable functions $f: X\to R$ with sup-norm. 

We also consider Banach spaces of bounded measures $\mu : \Sigma\to R$, with the norm equal to the total variation of the measure $\mu$ (but you can also use the equivalent sup-norm): 

$ba(X,\Sigma)$ is the space of finitely additive measures, 

$ca(X,\Sigma)$ is the space of countably additive measures. 

If $\mu\ge{0}$, then norm $||\mu||=\mu(X)$.

\begin{definition} \label{D1} 
	(\cite{YoHew1}). A finitely additive measure $\mu$, $\mu \ge 0$, is called purely finitely additive (pure charge, pure mean) if any countably additive measure $\lambda$ satisfying the condition $0\le\lambda\le\mu$ is identically zero. An alternating measure $\mu$ is called purely finitely additive if both components of its Jordan decomposition $\mu = \mu^{+}-\mu^{-}$ are purely finitely additive.
\end{definition}

If the measure $ \mu $ is purely finitely additive, then it is easy to see that it is equal to zero on every one-point set: $ \mu (\{x \}) = 0, \forall x \in X $. The converse, generally speaking, is not true, for example, for the Lebesgue measure on the segment $ [0, 1] $.

\begin {comment}
If the measure $ \mu $ is identically zero, then it can formally be considered both countably additive and purely finitely additive.

\end {comment}
 
\begin{theorem} 

(Alexandroff-Yosida-Hewitt decomposition, see \cite{YoHew1}).
Any finitely additive measure $ \mu $ can be uniquely decomposed into the sum $ \mu = \mu_{ca} + \mu_{pfa} $, where $ \mu_{ca} $ is countably additive and $ \mu_{pfa} $ is a purely finitely additive measure.

\end{theorem}

 Bounded purely finitely additive measures also form a Banach space $ pfa(X,\Sigma) $ with the same norm and $ ba(X,\Sigma) = ca(X,\Sigma) \oplus pfa(X,\Sigma) $ .
 
We denote the sets of nonnegative measures:

$V_{ba}=\{\mu\in{ba(X,\Sigma)}:\mu(X)\le 1\},$

$V_{ca}=\{\mu\in{ca(X,\Sigma)}:\mu(X)\le 1\},$

$V_{pfa}=\{\mu\in{pfa(X,\Sigma)}:\mu(X)\le 1\}.$
 
 Measures from these sets will be called probabilistic if $ \mu(X) = 1 $.
 
We also denote by $ S_{ba} $, $ S_{ca} $ and $ S_{pfa} $ the sets of all probability measures in $ V_{ba} $, $ V_{ca} $ and in $ V_{pfa} $, respectively.

\begin{definition}  \label{D2} 
The classical Markov chains (MC) on a measurable space $(X,\Sigma)$ are given by their transition function (probability kernel) $P(x,E), x\in X, E\in\Sigma$, under the usual conditions: 

\begin{enumerate}
	\item $0\le P(x,E) \le{1}, \forall{x}\in{X}, \forall{E}\in\Sigma$;
	\item $P(\cdot,E)\in{B(X,\Sigma)}, \forall{E}\in\Sigma$;
	\item $P(x,\cdot)\in{ca(X,\Sigma)}, \forall{x}\in{X}$;
	\item $P(x,X)=1, \forall{x}\in{X}.$
\end{enumerate}

\end{definition}

The numerical value of the function $ P (x, E) $ is the probability that the system will move from the point $ x \in X $ to the set $ E \in \Sigma $ in one step (per unit of time).

We emphasize that the transition function of the classical Markov chain is a countably additive measure in the second argument. 

We will also call such transition functions countably additive kernels.\\

The transition function generates two Markov linear bounded positive integral operators: 

$T: B(X,\Sigma)\to{B(X,\Sigma)}, (Tf)(x)=Tf(x)=\int\limits_X f(y)P(x,dy),$

$\forall{f\in{B(X,\Sigma)}},\forall{x}\in{X};$ \medskip

$A: ca(X,\Sigma)\to{ca(X,\Sigma)}, (A\mu)(E)= A\mu(E)=\int\limits_X P(x,E)\mu(dx),$

$\forall{\mu\in{ca(X,\Sigma)}},\forall{E}\in\Sigma.$ \medskip

Let the initial measure be $\mu_0\in{S_{ca}}$. Then the iterative sequence of countably additive probability measures $\mu_{n}=A\mu_{n-1}\in{S_{ca}},n\in{N}$, is usually identified with the Markov chain. 
We will call $ \{\mu_ {n}\}$ {\it a Markov sequence of measures}.\\
 
 Topologically conjugate to the space $ B (X, \Sigma) $ is (isomorphically) the space of finitely additive measures: $ B^* (X, \Sigma) = ba (X, \Sigma) $ (see, for example, \cite {DS1}). In this case, the operator $ T^* \colon ba (X, \Sigma) \to {ba (X, \Sigma)} $ serves as topologically conjugate to the operator T, which is uniquely determined by the well-known rule in terms of integral `` scalar products '':    
$$
\langle T^*\mu, f \rangle=\langle \mu, Tf\rangle, \forall{f\in{B(X,\Sigma)}},\forall{\mu\in{ba(X,\Sigma)}}.
$$ 

The operator $T^*$ is the only bounded continuation of the operator $A$ to the entire space $ba(X,\Sigma)$ while preserving its analytic form
$$
T^*\mu(E)=\int_X P(x,E)\mu(dx), \forall{\mu\in{ba(X,\Sigma)}}, \forall E\in\Sigma.
$$ 

The operator $T^*$ has its own invariant subspace $ca(X,\Sigma)$, i.e. $T^*[ca(X,\Sigma)] \subset{ca(X,\Sigma)}$, on which it coincides with the original operator $A$. The construction of the Markov operators $T$ and $T^*$ is now functionally closed. We shall continue to denote the operator $T^*$ as $A$.

In such a setting, it is natural to admit to consideration also the Markov sequences of probabilistic finitely additive measures
$$
\mu_0\in{S_{ba}}, \mu_n = A\mu_{n-1}\in{S_{ba}}, n\in{N},
$$
keeping the countable additivity of the transition function $P(x,\cdot)$ with respect to the second argument.

Despite this circumstance, the image $ A\mu $ of a purely finitely additive measure $ \mu $ can remain purely finitely additive, i. e., generally speaking,

$$A[ba(X,\Sigma)] \not \subset{ca(X,\Sigma)}.$$ 

{\bf Reference 2.1.} The integral over a finitely additive measure, usually called the Radon integral, is constructed according to the same scheme as the Lebesgue integral over the Lebesgue measure. Its construction is developed in \cite {YoHew1} and also, in a more modern form, in \cite {RaoRao1}. Note that if the original space $ X $ is countable and the measure $ \mu $ is not countably additive, then the integral on $ X $ cannot be replaced by a sum (series). Such integrals have other features as well. 

\begin{definition}  \label{D2} 
	If $A\mu = \mu$ holds for some positive finitely additive measure $\mu$, then we call such a measure {\it invariant} for the operator $A$ (and for the Markov chain).
	
	An invariant probability countably additive measure is often called the {\it stationary distribution} of a Markov chain.
\end{definition}

The question of the existence of invariant measures and their properties is one of the main questions in the theory of Markov chains.

We denote the sets of all nonzero invariant measures for the  operator $A$:

$\Delta_{ba}=\{\mu\in{V_{ba}}: \mu=A\mu\}$,

$\Delta_{ca}=\{\mu\in{V_{ca}}: \mu=A\mu\}$, 

$\Delta_{pfa}=\{\mu\in{V_{pfa}}: \mu=A\mu\}$.

The classical Markov chain with a countably additive transition probability may or may not have invariant countably additive probability measures, i.e. possibly $ \Delta_ {ca} = \emptyset $ (for example, for a symmetric walk on $ Z $).\\

In \cite [theorem 2.2] {Si1}
   Z. Shidak proved that any countably additive MC on an arbitrary measurable space $ (X, \Sigma) $ with an operator extended to the space of finitely additive measures has at least one invariant finitely additive measure, i.e., always $ \Delta_ {ba} \ne \emptyset $.

Z. Shidak in \cite [theorem 2.5] {Si1}
 in the general case also established for such MC that if a finitely additive measure $ \mu $ is invariant, $ A \mu = \mu $ and $ \mu = \mu_{ca} + \mu_{pfa} $ is its decomposition into countably additive and purely finitely additive components, then each of them is also invariant: $ A \mu_{ca} = \mu_{ca} $, $ A \mu_{pfa} = \mu$

 As already noted in the Introduction, in 1981 S. Ramakrishnan in his work \cite {Rama1} introduced the concept of finitely additive Markov chains generated by finitely additive (in the second argument) transition probabilities (kernels).
 
\begin{definition}  \label{D2.4}  A transition function of a finitely additive MC on an arbitrary (phase) measurable space $ (X, \Sigma) $ is a function $ P (x, E), x \in X, E \in \Sigma $, for which the conditions (1), (2) and (4) from definition 2.2, and instead of condition (3) condition $(3^{\prime})$ Is satisfied: $ P (x, \cdot) \in {ba (X, \Sigma)}, \forall {x} \in {X} $. We will also call such transition functions finitely additive.
\end{definition}

 The finitely additive transition function $ P (x, E) $ also generates two integral operators:
operator
$ T: B (X, \Sigma) \to {B (X, \Sigma)} $
  and operator
   $ A: ba (X, \Sigma) \to {ba (X, \Sigma)} $,
   with the same analytical form, with $ T ^ {*} = A $. However, in this case, generally speaking, the operator $ A $ does not transform countably additive measures into the same ones, that is, $ Aca (X, \Sigma) \not \subset ca (X, \Sigma) $.
   
  It is natural to consider the decomposition of such transition functions (kernels) into two components - countably additive and purely finitely additive. Such components may no longer have a direct probabilistic meaning.

To define such kernels, we will take as a basis Definition 2.2 and some information from Revuse's book \cite [Chapter 1, $ \S $ 1] {Rev1} and transfer them to the finitely additive case.

	\begin {definition}
A numerical function $ P (x, E) $ of two variables $ x \in X $ and $ E \in \Sigma $ is called a sub-Markov countably additive kernel if conditions (1), (2), (3) from Definition 2.2 are satisfied .
		
Such a kernel is called a Markov countably additive kernel if condition (4) from the same Definition 2.2 is additionally satisfied.
		
	\end{definition} 

Similarly, we introduce the terms sub-Markov and Markov kernels for the cases when the kernel $ P (x, \cdot) $ is finitely additive or purely finitely additive in the second argument for each $ x \in X $.

We can say that in this case we replace condition (3) in Definitions 2.2 and 2.4 by the conditions

($3^\prime$) $P(x, \cdot) \in ba(X, \Sigma), \forall x \in X$, and on the condition

($3^{\prime\prime}$) $P(x, \cdot) \in pfa(X, \Sigma), \forall x \in X$, respectively.

The integral operators $ T $ and $ A $ in spaces of functions and measures generated by a sub-Markov (Markov) kernel will also be called sub-Markov (Markov). \medskip

The already cited Alexandroff-Yosida-Hewitt theorem 2.1 (\cite {YoHew1}) on the decomposition of a finitely additive measure implies the following statement.
	
	\begin {proposition}
Let $ X $ be an infinite set and an arbitrary sigma-algebra of its subsets $ \Sigma $ contains all one-point sets. Any Markov finitely additive kernel $ P (x, E) $, on $ (X, \Sigma) $, is uniquely representable as the sum of its countably additive and purely finitely additive components:
	$$
	P(x, E)=P_{ca}(x, E)+P_{pfa}(x, E),
	$$
where 
	 $P_{ca}(x, \cdot) \in ca(X, \Sigma), P_{pfa}(x, \cdot)\in pfa(X, \Sigma),$ for all   $x \in X$, $E\in \Sigma$.
\end{proposition} 

We cannot yet call the components $ P_ {ca} (x, \cdot) $ and $ P_ {pfa} (x, \cdot) $ sub-Markov kernels, since the $\Sigma$ -measurability of the functions $ P_ {ca} (\cdot, E) $ and $ P_ {pfa} (\cdot, E) $ for different $ E \in \Sigma $ and for an arbitrary sigma-algebra $ \Sigma $ is not guaranteed. 
Moreover, the original Markov kernel $ P (\cdot, E) $  is $\Sigma $ -measurable for any $ E \in \Sigma $ by definition.

If the components $ P_ {ca} (\cdot, E) $ or $ P_ {pfa} (\cdot, E) $ are immeasurable, then there are no sub-Markov operators $ T $ and $ A $ integrally expressed in terms of them.\\

{\bf Reference 2.2.}
The question of measurability with respect to the first argument of two components in the decompositions of the Markov kernel in proposition 2.1 was pointed out by one of the authors of this article in the paper \cite {Zhd01}. It was hypothesized that immeasurable decompositions exist. This problem was solved by A.E. Gutman and A.I. Sotnikov in their work \cite {GuSot1}.
 
  They proved a number of theorems on the singularities of the decompositions of transition functions (kernels) into the sum of their countably additive and purely finitely additive components in different cases. It was proved that non-measurable decompositions exist, in particular, on the segment $ [0,1] $ with Lebesgue sigma-algebra.
 
Later, A.I. Sotnikov in his work \cite {Sot2} constructed a class of strongly additive transition functions in which both of their decompositions components are measurable.

In this paper, we will use another possibility of ensuring the measurability of the components in the decompositions of the finitely additive Markov kernel. This is an introduction in the next subsection of discrete topology in an arbitrary MC phase space.

\section {\bf  Markov kernels in discrete space} 

      In the theory of Markov chains the term "discrete" is used in different senses, and is applied to both the time parameter and the state space of the MC. We will use the classical definition from functional analysis (see, for example, \cite [Ch. II, $ \S $ 5, item 1, example 2] {KolFom}), which is also used in some papers on the theory of MC.

\begin{definition}
		A topological space $ (X, \tau) $ is called discrete if all its subsets are simultaneously open and closed (clopen), that is, the topology $ \tau = 2 ^ X $ - the set of all subsets of the set $ X $.
\end{definition} 

Such a topology in $ X $ is generated by the discrete metric $ d (x, y) $ equal to $ 1 $ for $ x \ne y $ and equal to $ 0 $ for $ x = y $.

In discrete space, all points are metrically isolated.

Discrete metric (and topology) can be introduced in any set $ X $. In particular, the discrete topology can be introduced in all "principal" number sets: $ N, Z, Q$, $ [0,1] $, $ R = R ^ {1} $, as well as in $ R ^ {m} (m \in N) $, transforming them into discrete spaces.

If a topological space is discrete, then, obviously, its Borel sigma-algebra $ \mathfrak {B} = \tau = 2 ^ X $. This sigma-algebra contains all subsets of the set $ X $. Such a sigma-algebra in $ X $ is also called discrete. We will denote it by $ \Sigma_ {d} $.

S. Ramakrishnan's work \cite {Rama1} uses a similar definition of the discrete phase space of a MC.

If the space $ X $ is discrete, then, obviously, any bounded numerical function $ f: X \to {R} $ is measurable with respect to the discrete sigma-algebra $ \Sigma_ {d} $, that is, $ f \in B (X , \Sigma_ {d}) $. In particular,  are  $ \Sigma_ {d} $ - measurable in the first argument and the components $ P_ {ca} (\cdot, E) $ and $ P_ {pfa} (\cdot, E) $ of the CM transition function in Proposition 2.1. for all $ E \in \Sigma $.

      Note that all numeric functions $ f: X \to {R} $ on any discrete space $ (X, \Sigma_ {d}) $ are continuous in the discrete topology $ \tau = 2 ^ X $.

Let us introduce the concept of a measure atom, known in different versions (we just need to use a simplified version of its definition).

\begin{definition}  
	Let $ (X, \Sigma) $ be an arbitrary measurable space and $ \mu: \Sigma \to R $ some countably additive measure. An element $ x \in X $ is called an atom of the measure $ \mu $ if $ \mu (\{x \})\ne 0 $. If a bounded measure $ \mu, \mu \ge 0, $ has a support (set of full measure) $ D \in \Sigma $, consisting of a finite or countable family of its atoms, then such a measure is called atomic (discrete). Moreover, $ D = \{x_ {1}, x_ {2}, \dots \} $ and $ \mu (D) = \sum_ {n}\mu (\{x_ {n} \}) = \mu (X) $.

\end{definition} 

The atomic measure $ \mu \ge 0 $ can be represented as follows
$$
\mu(E)=\sum_{n} \alpha_{n} \delta_{x_{n}}(E),
$$
where $ E \in \Sigma $, $ \delta_ {x_ {n}} $ are Dirac measures concentrated at the points $ x_ {n} $, and $ \sum_ {n} \alpha_ {n} = \mu (D ) = \mu (X) $.

Note that a countably additive measure on a nondiscrete measurable space $ (X, \Sigma) $ may not have atoms, for example, the Lebesgue measure on $ ([0,1], \mathfrak {B}) $.
It is also true that any purely finitely additive measure on any measurable space has no atoms.

If the set $ X $ is countable and $ \Sigma = \Sigma_ {d} $, then, obviously, any bounded  countably additive measure $ \mu $ on $ (X, \Sigma_ {d}) $ is atomic.\\

Now we want to find out how countably additive measures are arranged on an arbitrary discrete space $ (X, \Sigma_ {d}) $. This question in a broader setting is considered, for example, in N. Bourbaki's book \cite [Ch. III, paragraphs 1-2] {Burb1}. A locally compact topological space is taken as the initial space $ X $. Countably additive measures are defined as linear continuous functionals on the space of continuous functions. Definitions of a discrete space, a discrete (atomic) measure and its support are given, which differ from those given above.
After proving a number of propositions (theorems), in \cite [Ch. III, paragraph 2, item 5] {Burb1}, the following statement is formulated: `` on a discrete space, any measure is discrete '' (here are meant countably-additive measures).

To apply this statement in this work, we would need to give precise definitions of the above and other concepts and translate them into our language. Therefore, in theorem 3.2 below, we give our proof of the above statement from \cite {Burb1} in our definitions and refine it.

However, for this we need one well-known and nontrivial theorem of Ulam, stated, for example, in \cite [Ch. 5, theorems 5.6 and 5.7] {O1} and, in more detail, in \cite [Volume 1, theorem 1.12.40, and corollary 1.12.41)] {Boga1}. 
We present this theorem under the condition that the continuum hypothesis is accepted, i.e. we assume that $ \aleph_ {1} = c $ (continuum).

\begin{theorem} 
(Ulam).
A finite countably additive measure $ \mu $ defined on all subsets of the set $ X $ of cardinality $ \aleph_ {1} $ ($ c $, continuum) is identically zero if it is zero for each one-point subset.
\end{theorem} 

 Obviously, Ulam's theorem trivially holds for sets $ X $ with countable cardinality $ \aleph_ {0} $.
 
 We will continue to assume that the continuum hypothesis is true.

\begin{comment} 
In the books \cite {O1} and \cite {Boga1} it is noted that the (extended) Ulam theorem is true and for higher, so-called `` immeasurable '' cardinalities of the set $ X $. Immeasurable cardinality includes all cardinality from an ordered cardinality scale: $ \aleph_ {0} $, $ \aleph_ {1} = c $, $ \aleph_ {2} $, etc. There is still no example of a set with `` measurable '' power.
\end{comment}

 \begin{definition} 
  A measurable space $ (X, \Sigma_ {d}) $ will be called an arbitrary discrete space if $ \Sigma_ {d} = 2 ^ {X} $, and the set $ X $ has an arbitrary immeasurable cardinality (including from the ordered cardinality scale ). In other words, we consider only those discrete spaces for which the (extended) Ulam theorem is valid.
 \end{definition} 

 Now let us prove the following promised theorem.

\begin{theorem} 
	  Any non-zero non-negative bounded countably additive measure $ \mu: \Sigma_ {d} \to R $, on an arbitrary discrete space $ (X, \Sigma_ {d}) $ is atomic (discrete), and has a finite or countable support $ D = \{x_ {1}, x_ {2}, \dots\} \subset X $, for which $ \mu (\{x_ {n}\}) = \alpha_ {n}> 0, n \in N $, $ \sum_ {n} \alpha_ {n} = \mu (D) = \mu (X) $, $ \mu (X \setminus D) = 0 $.
\end{theorem} 

Proof. For a countable set $ X $, the assertions of the theorem are trivially fulfilled.
 
 Now let the set $ X $ have uncountable cardinality. Consider an arbitrary bounded nonnegative countably additive measure $ \mu: \Sigma_ {d} \to R $ for which $ 0 <\mu (X) = \gamma <\infty $. Since the measure $ \mu $ is not identically zero, then, by theorem 3.1, the measure $ \mu $ has at least one one-point atom $ x_ {0} \in X $ such that $ \mu (\{x_ { 0} \}) = \alpha_ {0}> 0 $. We denote by $ D $ the set of all atoms of measure $ \mu $. As we have shown above, $ x_ {0} \in D $ and $ D \ne \O $. Let us prove that the set $ D $ is finite or countable.
 
We split the interval $ (0, \gamma] $ of possible nonzero values of the measure $ \mu $ into a countable family of disjoint intervals
$$
 (0, \gamma]= \cup_{n=1}^{\infty} (\frac{\gamma}{n+1}, \frac{\gamma}{n} ]. 
$$
 
We denote the inverse images of these intervals
$$
D_{n}=\{x\in X: \mu(\{x\})\in (\frac{\gamma}{n+1}, \frac{\gamma}{n} ] \}, D_{n}\in\Sigma_{d}, n=1, 2, \dots . 
$$

Then the sets $ D_ {n} $ are also pairwise disjoint and $ D = \cup_ {n = 1} ^ {\infty} D_ {n} $. Therefore, since the measure $ \mu $ is countably additive, then $ \mu (D) = \sum_ {n = 1} ^ {\infty} \mu (D_ {n}) $.

By construction, for any point $ x \in D_{n} $ performed $ \frac{\gamma}{n + 1} <\mu (\{x \}) \le \frac{\gamma}{n} $, $ n = 1, 2, \dots $. In addition, $ \mu (D_{n}) \le \gamma <\infty $ for all $ n = 1, 2, \dots $.

If any of the sets
$ D_ {n} $ was infinite, then, by virtue of the inequalities for $ \mu (\{x \}) $ for $ x \in D_ {n} $, it would be $ \mu (D_ {n}) = \infty. $ 
This contradiction implies that each set $ D_ {n} $ is finite or empty.

Therefore, the set $ D $, as a union of a countable (or finite) family of finite sets, is countable (or finite).

By construction, $ D \subset X $ and $ \mu (D) \le \mu (X) $. Let us prove that $ \mu (D) = \mu (X) $.
Since $ X $ has uncountable cardinality, and the set $ D $ is finite or countable, the set $ X \setminus D \ne \O $ and is also uncountable.

Suppose $ \mu(D)<\mu(X) $, i.e. $ \mu(X \setminus D)>0 $. By the hypothesis of the theorem, the set $ X $ is discrete. Consequently, the set $ X \setminus D $ is also discrete.

The restriction $ \mu_ {D} $ of the measure $ \mu $ from the set $ X $ to the set $ X \setminus D $ also satisfies all the requirements for the measure $ \mu $ under the conditions of the theorem. Since $ \mu (X \setminus D)> 0 $, then, again applying Ulam 3.1's theorem, we obtain that there exists a point $ y_ {0} \in X \setminus D $ such that $ \mu_ {D} ( \{y_ {0}\})> 0 $. But then $ \mu (\{y_ {0}\})> 0 $. Therefore, $ y_ {0} \in D $, and $ y_ {0} \ne x $ for any $ x \in D $. However, the set $ D $ was defined as the set of all points $ x \in X $ for which $ \mu (\{x \})> 0 $ is satisfied. We got a contradiction. Therefore, $ \mu (D) = \mu (X) $ and $ \mu (X \setminus D) = 0 $.

Since $ D $ is finite or countable, we renumber all its points and obtain the last statement of the theorem.

The theorem is proved.\\

{\bf Reference 3.1.} It may be surprising that all countably additive measures on the discrete segment $ [0,1] $ are only atomic. And where is the Lebesgue measure? Yes, it really does not exist on the discrete segment $ [0,1] $. In 1923 Stefan Banach proved that the Lebesgue measure defined on the Borel (generated by the Euclidean topology) sigma-algebra of the segment $ [0,1] $ cannot be countably additively extended to the sigma-algebra of all subsets of the segment $ [0,1 ] $. However, it can be extended to a finitely additive measure on a discrete sigma-algebra, and there are infinitely many such extensions. This issue is discussed in many sources, see, for example, \cite [volume 1, items 1.12.29 and 2.12.91] {Boga1}.

\begin{theorem}
 	  A finitely additive nonnegative measure $ \mu $ defined on an arbitrary discrete space $ (X, \Sigma_ {d}) $ is purely finitely additive if and only if the condition $ \mu (\{x \}) = 0 $ for all $ x \in X $ are fulfilled.
\end{theorem} 

 Proof. The necessity of the condition is obvious. Let us show its sufficiency. Let the condition be satisfied, but the measure $ \mu $ is not purely finitely additive. Then in its decomposition $ \mu = \mu_{ca} + \mu_{pfa} $ the countably additive component $ \mu_ {ca}\ne 0 $. In this case, by theorem 3.2, there exists a point $ x_ {1} \in X $ such that $ \mu_ {ca}(\{x_ {1}\})> 0 $ It follows from this contradiction that $ \mu $ is purely finitely additive. The theorem is proved.
 
Let us now return to Markov chains. Theorem 3.2 automatically implies the following statement.

\begin{theorem} 
	Let a countably additive sub-Markov kernel $ P (x, E) $ be given on an arbitrary discrete space $ (X, \Sigma_ {d}) $, and $ P (x, X)> 0 $ for all $ x \in X $ . Then for any $ x \in X $ the measure $ P (x, \cdot) $ is atomic and has a finite or countable support $ D (x) = \{x_ {1} (x), x_ {2} (x ), \dots \} $, for which\\
	$ P (x, \{x_ {n} (x) \}) = \alpha_ {n} (x)> 0, n \in N $,\\
	 and $ \sum_ { n} \alpha_ {n} (x) = P (x, D (x)) = P (x, X), P (x, X \setminus D (x)) = 0 $.
	\end{theorem} 
	
	\begin{corollary}
(from theorem 3.3). Let a finitely additive sub-Markov kernel $ P (x, E) $ be given on an arbitrary discrete space $ (X, \Sigma_ {d}) $. For any fixed $ x_ {0} \in X $ the measure $ P (x_ {0}, \cdot) $ is purely finitely additive if and only if $ P (x_ {0}, \{y \}) = 0 $ for all $ y \in X $ (including the case $ y = x_ {0} $).
\end{corollary}



{\bf Example  3.1.}
  Let $ X = [0,1]$  with Euclidean topology, $ \Sigma = \mathfrak {B} $, be a Borel sigma-algebra, and a countably additive MC given by the kernel $ P (x, E) = \lambda (E) $ for all $ x \in X $ and $ E \in \mathfrak {B} $, where $ \lambda $ is the Lebesgue measure. Such a MC corresponds to a sequence of independent uniformly distributed random variables on the segment $ [0,1] $. Obviously, it holds $ P (x, \{y \}) = 0 $ for all $ x, y \in X $. However, the phase space $ X = [0,1], \Sigma = \mathfrak {B} $ is not discrete, and theorem 3.3 is not applicable.\medskip
 
{\bf  Example 3.2.}
 Now we take on the same $ X $ the discrete sigma-algebra $ \Sigma_ {d} $. Consider a finitely additive MC defined by the kernel $ P (x, E) = \eta (E) $ for all $ x \in X $ and $ E \in \Sigma_ {d} $, where $ \eta $ is some purely finitely additive measure satisfying the conditions: $ \eta \ge 0 $, $ \eta (X) = 1 $ and $ \eta ((0, \varepsilon)) = 1 $ for all $ \varepsilon> {0} $. Then, obviously, the condition $ P (x, \{y \}) = 0 $ is also satisfied for all $ x, y \in X $, and theorem 3.3 is applicable.
 
The measure $ \eta $ in this example can be informally characterized as follows. It specifies a certain ``random variable'' that takes a value with probability 1, as close to point 0 as you like, but not at point 0.\medskip

We denote by $ P^{n} (x, E) $ the integral convolution of the kernel $ P^{1} (x, E) = P (x, E) $, $ n = 1,2,3,\ldots$.

\begin{corollary}
	 Let a sub-Markov purely finitely additive kernel $ P (x, E) $ be given on an arbitrary discrete space $ (X, \Sigma_ {d}) $. Then for all $ x, y \in X $ and $ n = 1,2,3, ... $, $ P^{n} (x, \{y\}) = 0 $ (including the case $ x = y $).
\end{corollary}

This statement is easily proved by induction.

  Generally speaking, the converse is not true. Here is a counter-example.\\
  
{\bf  Example  3.3.} 
Let some purely finitely additive probability measure $  \eta $ be given on the discrete space $(X, \Sigma_{d})$, where $ X = [0,1] $. Consider on $ ([0,1], \Sigma_ {d}) $ a finitely additive MC with the following rules for passing in one step:

   $P(0,\{1\})=1, P(x,E)=\eta(E)$ for all $x\in (0,1]$ and $E\in \Sigma_{d}$. In particular, $ P(x,\{y\})=\eta(\{y\}=0 $ for all $x\in (0,1]$ and $y\in [0,1]$.
   
    Performing the integral convolution of two kernels $ P (x, E) $, we obtain that $ P^{2} (x, \{y\}) = 0 $ for all $ x \in [0,1] $ and $ y \in [0,1] $. Moreover, $ P (0, \{1\}) = 1> 0 $. \medskip
   

 
 
  As noted above in Section 2, the operator $ A $ generated by the countably additive submarkov kernel transforms countably additive measures into the same ones, that is, $ A [ca (X, \Sigma)] \subset ca (X, \Sigma) $. This property is preserved for the particular discrete case $ \Sigma = 2^{X} $. However, if the measure $ \mu \in V_ {pfa} $ is purely finitely additive, then both cases are possible: $ A \mu \in V_ {ca} $ and $ A \mu \in V_ {pfa} $. However, the situation is different with a purely finitely additive kernel.

\begin{theorem} 
	 Let a purely finitely additive sub-Markov kernel $ P (x, E) $ be given on an arbitrary discrete space $ (X, \Sigma_ {d}) $. Then the sub-Markov operator $ A $ generated by this kernel transforms all finitely additive measures into purely finitely additive measures, that is, $ A [ba (X, \Sigma_ {d})] \subset pfa (X, \Sigma_ {d} ) $, in particular, $ A [ca (X, \Sigma_ {d})] \subset pfa (X, \Sigma_ {d}) $ and $ A [pfa (X, \Sigma_ {d})] \subset pfa (X, \Sigma_ {d}) $.
\end{theorem} 

Proof.
Let the finitely additive measure $ \mu \in V_ {ba} $ and $ \mu (X)> 0 $. We denote the measure by $ \eta = A \mu $. It is clear that the measure $ \eta $ is also finitely additive.

If $ \eta (X) = 0 $, that is, $ \eta \equiv 0 $ (this is possible), then it can be considered purely finitely additive (see remark 2.1) and the theorem is true.

Let $ \eta (X)> 0 $. Take its decomposition $ \eta = \eta_ {ca} + \eta_ {pfa} $ into a countably additive component $ \eta_ {ca} $ and a purely finitely additive component $ \eta_ {pfa} $.

If $ \eta_ {ca} (X) = 0 $, then the measure is $ \eta = \eta_ {pfa} \in V_ {pfa} $ and the theorem is proved.

Suppose that the countably additive measure $ \eta_ {ca} (X)> 0 $. Then, by theorem 3.2, the measure $ \eta_ {ca} $ has at least one atom $ a \in X $, $ \eta_ {ca} (\{a \}) = \gamma> 0 $.
Since the measure $ \eta_ {pfa} $ is purely finitely additive, then $ \eta_ {pfa} (\{a\}) = 0 $.

By the hypothesis of the theorem, all kernels $ P (x, \cdot) $ are purely finitely additive for all $ x \in X $. Such measures vanish on any one-point set. Therefore, $ P (x, \{a\}) = 0 $ for all $ x \in X $. Hence it follows that
$$
\gamma=\eta_{ca}(\{a\})=\eta_{ca}(\{a\})+ 0 =\eta_{ca}(\{a\})+\eta_{pfa}(\{a\})=\eta(\{a\})= A\mu(\{a\})
$$
$$
=\int_{X}P(x, \{a\})\mu(dx)=\int_{X}0 \cdot \mu(dx)=0.
$$

We got a contradiction. Therefore, $ \eta_ {ca} (X) = 0 $, and the measure $ \eta = \eta_ {pfa} = A \mu $ is purely finitely additive. The theorem is proved.
    
     Preliminary version of theorem 3.5. in the particular case $ X = Z $ was proved in the work of one of the authors of this article \cite {XAK1}.
    

Now, using the discrete topology in $ X $, we can complete proposition 2.1.

\begin{proposition}  
Let an arbitrary discrete space $ (X, \Sigma_ {d}) $ be given. Any Markov finitely additive kernel $ P (x, E) $ on $ (X, \Sigma_ {d}) $ is uniquely representable as the sum of a sub-Markov countably additive kernel $ P_{ca}(x, E) $ and a sub-Markov purely finitely additive kernel $ P_{pfa}(x, E) $:
$$
P(x, E)=P_{ca}(x, E)+P_{pfa}(x, E),
$$
where

 $P_{ca}(x, \cdot) \in ca(X, \Sigma_{d}), P_{pfa}(x, \cdot)\in pfa(X, \Sigma_{d}),$  
 and $P_{ca}(\cdot, E) \in B(X, \Sigma_{d})$, 
 $ P_{pfa}(\cdot, E)\in B(X, \Sigma_{d})$ for all $x \in X$ and $E\in \Sigma_{d}$.
\end{proposition}

The last inclusions mean that the kernels $ P_ {ca} (\cdot, E) $ and $ P_ {pfa} (\cdot, E) $  are $ \Sigma_ {d} $ - measurable in the first argument for all $ E\in \Sigma_ {d} $.

Proposition 3.1 makes it possible to introduce integral sub-Markov operators generated by the corresponding measurable subkernels.

Let's denote these operators $ A_{ca} $ and $ A_{pfa} $:
$$
\eta_{ca}(E)=A_{ca}\mu(E)=\int_{X}P_{ca}(x, E)\mu(dx), E\in \Sigma_{d}, \eta_{ca}=A_{ca}\mu,
$$
$$
\eta_{pfa}(E)=A_{pfa}\mu(E)=\int_{X}P_{pfa}(x, E)\mu(dx), E\in \Sigma_{d}, \eta_{pfa}=A_{pfa}\mu,
$$
where
$\mu, \eta_{ca}, \eta_{pfa} \in ba(X, \Sigma_{d}), A=A_{ca}+A_{pfa}$.



Obviously, both sub-Markov operators $ A_{ca} $ and $ A_{pfa} $ are linear, bounded (continuous), $\|A_{ca} \| \le 1 $, $ \|A_{pfa} \| \le 1 $, and are positive.

As we have already found out,

 $$A_{ca}[ca(X, \Sigma_{d})] \subset ca(X, \Sigma_{d}), A_{pfa}[ba(X, \Sigma_{d})] \subset pfa(X, \Sigma_{d}).$$ 
 




\begin{corollary} 
The following inclusions are true for superpositions of operators $A_{ca}$ and $A_{pfa}$:
    \begin{enumerate}
      \item $A_{ca} \cdot A_{ca}[ca(X, \Sigma_{d})] \subset ca(X, \Sigma_{d});$ 
        \item $A_{pfa} \cdot A_{pfa}[ba(X, \Sigma_{d})] \subset pfa(X, \Sigma_{d});$ 
       \item $A_{pfa} \cdot A_{ca} [ba(X, \Sigma_{d})] \subset pfa(X, \Sigma_{d}).$
    \end{enumerate}
\end{corollary}

\begin{comment}  
The operators $ A_{ca} $ and $ A_{pfa} $, generally speaking, are non-commutative, i.e.  $A_{ca} \cdot A_{pfa}$ $\neq$  $A_{pfa} \cdot A_{ca}$.
\end{comment}

 
 
 
\newpage



 \section {\bf Invariant measures of Markov kernels on discrete space}

  In the paper \cite [Chapter I, \S 5, theorem 5.3] {Zhd01} the following statement was proved.

\begin{theorem}
(Zhdanok) For any Markov chain with a Markov finitely additive kernel $ P (x, E) $ on an arbitrary measurable space $ (X, \Sigma) $ there exists an invariant probability finitely additive measure $ \mu = A \mu \in S_ {ba} $, that is, $ \Delta_ {ba} \ne \O $.

\end{theorem} 

Earlier, a similar theorem (in the language of strategies) was proved by Ramakrishnan \cite [p. 8, theorem 2] {Rama1}, but in the special case of a discrete phase space. In our theorem 4.1 given above, no restrictions on the phase space are assumed.

Now let on an arbitrary discrete space $ (X, \Sigma_ {d}) $, $ \Sigma_ {d} = 2^X $, a Markov chain with a Markov finitely additive kernel $ P (x, E) $ is given. We have previously identified two special "extreme" cases. The first is when the kernel $ P (x, \cdot) $ is a countably additive measure for every $ x \in X $. The second is when the kernel $ P (x, \cdot) $ is a purely finitely additive measure for all $ x \in X $.

The first case has already been considered in the previous paragraphs of this article and studied in a number of works by various authors.

Consider now the second special case.

\begin{theorem} 
Let a Markov chain with a purely finitely additive kernel $ P (x, E) $ be given on an arbitrary discrete space $ (X, \Sigma_ {d}) $. Then, for the Markov operator $ A $ generated by it, there exists an invariant probabilistic finitely additive measure $ \mu = A \mu \in V_ {ba} $ and all its invariant measures are purely finitely additive, that is, $ \Delta_ {ba} = \Delta_ {pfa} \ne \O $ and  $\Delta_{ca}= \O$.
\end{theorem} 

 Proof. Theorem 4.1 is proved for any sigma-algebra $ \Sigma $ subsets of $ X $ and for any Markov finitely additive kernel. Hence, it is also true for the discrete sigma-algebra $ \Sigma_ {d} = 2^X $ and for a purely finitely additive kernel.
  
Therefore, under the conditions of the present theorem, for the operator $ A $ there exists an invariant probabilistic finitely additive measure $ \mu = A \mu $, defined on the discrete space $ (X, \Sigma_ {d}) $.

It follows from theorem 3.5 that the measure $ \mu $ and all other invariant measures of the operator $ A $ are purely finitely additive. The theorem is proved.


Obviously, in order to further study the properties of finitely additive MCs, taking into account the form of their decomposition, it is necessary to single out some natural non-trivial special cases. In this paper, we will consider one such option.

\begin{definition} 
	We call a finitely additive MC on an arbitrary discrete space $ (X, \Sigma_ {d}) $ combined if its transition function in the decomposition 
	$$P(x,E)=P_{ca}(x,E)+P_{pfa}(x,E),$$
	$$ P_{ca}(x,\cdot) \in ca(X, \Sigma_{d}), P_{pfa}(x,\cdot) \in pfa(X, \Sigma_{d}),$$ 
	satisfies the conditions:
	$$P_{ca}(x,X)=q_{1},  P_{pfa}(x,X)=q_{2} \text{ for all } x \in X,$$
	where $0\le q_{1}, q_{2}\le1$,  $q_{1}+q_{2}=1$.
\end{definition} 

Let the finitely additive MC be combined. Then, as shown in the comments to proposition 3.1, its Markov operator $ A $ can also be represented as the sum $ A = A_{ca} + A_{pfa} $ of its two components generated by the sub-Markov kernels $ P_{ca}(x, E ) $ and $ P_{pfa}(x, E) $. Wherein $\|A_{ca}\|= q_{1}, \|A_{pfa}\|= q_{2}$.

\begin{definition} 
	A combined MC is called non-degenerate if its decomposition from definition 4.1 holds $ 0 <q_ {1}, q_ {2} <1 $, and degenerate if $ q_ {1} = 0 $ or $ q_ {2} = 0 $.
\end{definition} 

Above, in Section 2 and in theorem 4.2, we describe the existence of invariant measures and their types for countably additive and purely finitely additive MCs. They, by definition 4.2, are degenerate cases of combined MCs.

Let the MC not be degenerate. Let's take functions 
$$\tilde{P}_{ca}(x,E)=\frac{1}{q_{1}}P_{ca}(x,E), \tilde{P}_{pfa}(x,E)=\frac{1}{q_{2}}P_{pfa}(x,E).$$

Then the functions $ \tilde{P}_{ca} (x, E) $ and $ \tilde{P}_{pfa} (x, E) $ satisfy definition 2.1 and are transition functions (Markov kernels) of the corresponding Markov operators
 $$\tilde{A}_{ca}=\frac{1}{q_{1}}A_{ca},  \tilde{A}_{pfa}=\frac{1}{q_{2}}A_{pfa}.$$
 
Therefore, the Markov operator $ A $ of the combined MC is a linear combination $ A = q_ {1} \tilde{A}_{ca} + q_ {2} \tilde{A}_{pfa} $ for two Markov operators $ \tilde{A}_{ca} $ and $ \tilde{A}_{pfa} $ (hence the name of such MCs and operators in definition 4.1 is taken).

Recall that, by theorem 4.1, any, including combined, finitely additive MC has an invariant finitely additive measure.

\begin{theorem} 
	The combined non-degenerate finitely additive  MC on an arbitrary discrete space $ (X, \Sigma_ {d}) $ has no nonzero invariant countably additive measures, that is, $ \Delta_ {ca} = \O $.
\end{theorem} 

We carry out the proof by contradiction.

Suppose that $ \mu = A \mu \in S_ {ca} $, i.e. the invariant measure $ \mu $ is countably additive). 
	
	Then $ \mu = A \mu = (A_{ca} + A_{pfa}) \mu = A_{ca} \mu + A_{pfa} \mu $, where $ A_{ca} $ is countably additive , and $ A_{pfa} $ are purely finitely additive components of the operator $ A $. Then $ A_{ca} \mu $ is also a countably additive measure, and $ A_{ca} \mu (X) = q_ {1}> 0 $, that is, the measure $ A_{ca} \mu $ is nonzero. By theorem 3.5, the measure$ A_{pfa} \mu $ is purely finitely additive and also nonzero, $ A_{pfa} \mu (X) = q_ {2}> 0 $.

Consequently, the measure $ \mu $ has a nonzero purely finitely additive component $ A_{pfa} \mu $, and is not countably additive. The resulting contradiction proves the theorem.
	
\begin{theorem} 
	Let on an arbitrary discrete space $ (X, \Sigma_ {d}) $ a combined non-degenerate finitely additive MC with invariant probability finitely additive measure $ \mu = A \mu \in S_ {ba} $ is given. Let $ \mu = \mu_ {ca}+ \mu_{pfa}$ be its decomposition into countably additive $ \mu_ {ca}$ and purely finitely additive $ \mu_{pfa}$ components, and $ \mu_ {1 } \ne 0 $, $ \mu_{pfa}\ne 0 $.
	
	Then the measures $ \mu_ {ca}$ and $ \mu_{pfa}$ are not invariant for the operator $ A $, that is, $ \mu_ {ca}\ne A \mu_ {ca}$ and $ \mu_{pfa} \ne A \mu_{pfa}$.
\end{theorem} 

Proof.	
If $ \mu_ {ca}= A \mu_ {ca}$, then the MC has an invariant countably additive measure $ \mu_ {ca}$, which contradicts theorem 4.3. Therefore, $ \mu_ {ca}\ne A \mu_ {ca}$.

	From the decomposition of the measure $ \mu $, we have:
	$$\mu_{ca}+\mu_{pfa}=\mu=A\mu=A\mu_{ca}+A\mu_{pfa}.$$ 	
	
	Hence, if $ \mu_{pfa}= A \mu_{pfa}$, then $ \mu_ {ca}= A \mu_ {ca}$. And this contradicts the previous statement. Therefore, $ \mu_{pfa}\ne A \mu_{pfa}$.
		
	The theorem is proved.\medskip
		
	The initial versions of theorems 4.3 and 4.4 in the particular case $ X = Z $ were presented in the work of one of the authors of this article \cite {XAK1}. Combined MCs are called mixed there.

Recall that by Shidak's theorem \cite [theorem 2.5] {Si1}, for a MC with a countably additive kernel in a similar decomposition of the invariant measure $ \mu = \mu_ {ca}+ \mu_{pfa}$, performed $ \mu_ {ca} = A \mu_ {ca}$ and $ \mu_{pfa}= A \mu_{pfa}$. The difference between such CMs and combined ones turned out to be very significant.

Let us give an example to illustrate the last two theorems.\\

{\bf Example  4.1.} 
Consider on the segment $ X = [0,1] $ with discrete sigma-algebra $ \Sigma_ {d} $ a combined finitely additive MC with kernel
 $$P(x,E)=P_{ca}(x,E)+P_{pfa}(x,E),$$
where $ P_{ca}(x, E) $ is its countably additive component, and $ P_{pfa}(x, E) $ is purely finitely additive. These components are set according to the rules:

$ P_{ca}(x, E) = \frac {1} {2} \delta_ {0} (E) $ for all $ x \in X $ and $ E \subset X $, where $ \delta_ {0 } $ -- Dirac measure at point 0;

 $ P_{pfa}(x, E) = \frac {1} {2} \eta (E) $ for all $ x \in X $ and $ E \subset X $, where $ \eta $ is some fixed  purely finitely additive measure from $S_{pfa}$. 
 For clarity, we take the measure $ \eta $ from the family of purely finitely additive measures satisfying the condition $ \eta ((0,\varepsilon)) = 1 $ for any $ \varepsilon > 0 $.
 
 Moreover, $ P_{ca} (x, X) = \frac{1}{2} = q_{1} $ and $ P_{pfa}(x, X) = \frac{1}{2} = q_{2} $ for all $ x \in X $.

\begin{figure}[ht!] 
	\centering
	\includegraphics[width=95mm]{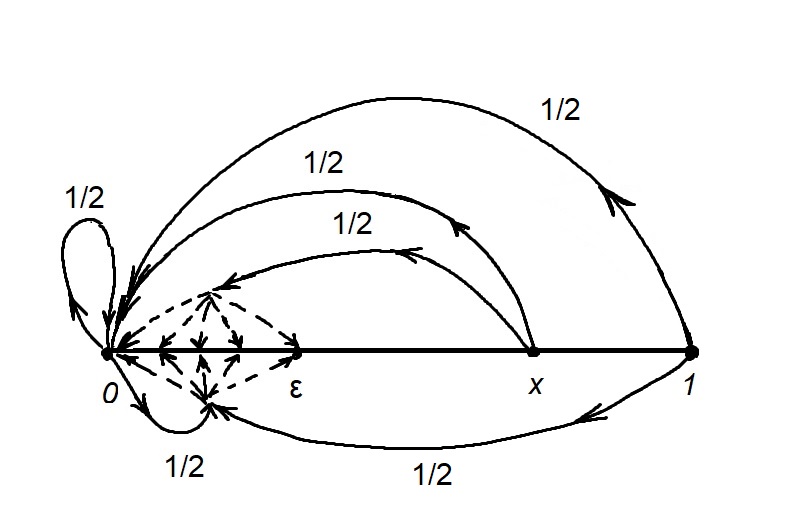}
	\caption{Phase portrait of the MC from Example 4.1.} 
	\label{ris:1}
\end{figure}

Essentially, all this means that a Markov chain in one step can move from any point $x\in X$ to point 0 with probability $ \frac{1}{2} $, and to any set $ E \subset X $ with probability $ \frac{1}{2}\eta (E) $. 
In particular, from any point $x\in X$, the system can move with probability $\frac{1}{2}$ to the open interval $(0,\varepsilon)$ for every $\varepsilon \in (0,1)$.
The phase portrait of such a MC with an arbitrary $\varepsilon \in (0,1)$ is shown in the Figure~\ref{ris:1}.

Take an arbitrary (initial) finitely additive probability measure $ \mu \in S_{ba} $. Then for any $ E \subset X $ the following holds:
$$A\mu(E) = \int_{X}P(x,E)d\mu(x)
= \int_{X}P_{ca}(x,E)d\mu(x) + \int_{X}P_{pfa}(x,E)d\mu(x) 
$$
$$
= \frac{1}{2}\int_{X}\delta_{0}(E)d\mu(x) + \frac{1}{2}\int_{X}\eta(E)d\mu(x)
$$
$$
= \frac{1}{2}\delta_{0}(E)\cdot \mu(X) + \frac{1}{2}\eta(E) \cdot \mu(X) 
= \frac{1}{2}\delta_{0}(E) + \frac{1}{2}\eta(E).
$$

Hence,
 $A\mu =  \frac{1}{2}\delta_{0} + \frac{1}{2}\eta$ 
  for any initial measure $\mu$.
 
 If we put
 $\mu =\frac{1}{2}\delta_{0} + \frac{1}{2}\eta$, 
  we get  $A\mu =\mu$.
 
 Obviously, this is the only invariant probabilistic finitely additive measure for a given MC.
 
The measures $ \mu_ {ca}= \frac {1} {2} \delta_ {0} $ and $ \mu_{pfa}= \frac {1} {2} \eta $ are nonzero components of the measure $ \mu $ - countably additive and purely finitely additive, respectively, and $ \mu = \mu_ {ca}+ \mu_{pfa}$. Thus, theorem 4.3 is confirmed. It is also obvious that $ A \mu_ {ca}= \mu \ne \mu_ {ca}$ and $ A \mu_{pfa}= \mu \ne \mu_{pfa}$. Therefore, this example also confirms theorem 4.4.\\

In the combined non-degenerate decomposition $ A = A_{ca} + A_{pfa} $ of the finitely additive operator $ A $, its countably additive component $ A_{ca} $ and the purely finitely additive component $ A_{pfa} $ enter equally. One might suppose that theorem 4.3 would also be valid for a purely finitely additive invariant measure. However, it is not. Let us give a corresponding counterexample.\\

{\bf  Example 4.2.} 
We consider a finitely additive combined MC on a discrete segment $ [0,1] $ under the same conditions as in example 4.1, but with a different countably additive component of its kernel:

$P_{ca}(x,E)=\frac{1}{2} \delta_{x}(E)$ for all $x\in X$ and $E\subset X$, where  $\delta_{x}$ -- Dirac measure at point $x$.

\begin{figure}[ht!] 
	\centering
	\includegraphics[width=95mm]{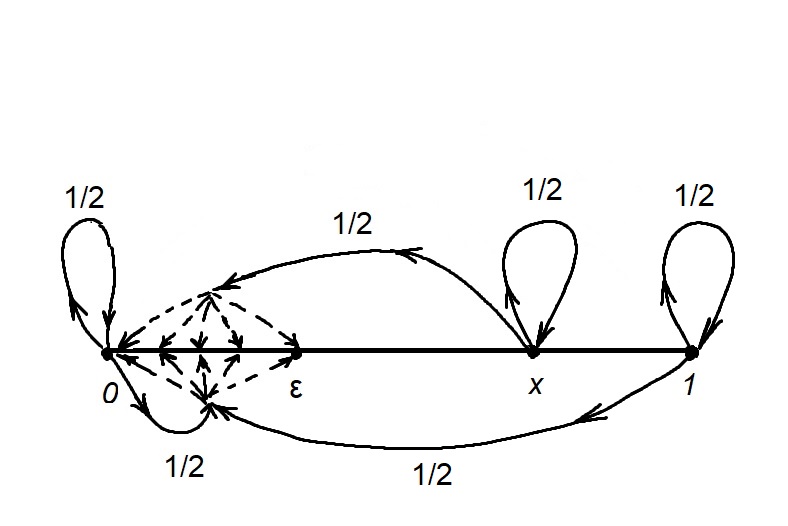}
	\caption{Phase portrait of the MC from Example 4.2.} 
    \label{ris:2}
\end{figure}

Meaningfully, this means that in one step the Markov system can go from any $ x \in X $ to the point $ x $, i.e. go into itself with probability $ \frac{1}{2} $, and into any set $ E \subset X $ with probability $ \frac{1}{2} \eta(E) $. 
In particular, the probability $P_{pfa}(x,(0,\varepsilon))=\frac{1}{2}$ for any  $\varepsilon \in (0,1)$.
The phase portrait of such a MC with an arbitrary $\varepsilon \in (0,1)$ is shown in the Figure~\ref{ris:2}.

Obviously, this MC is a combined non-degenerate chain.

 Let us perform integral transformations for an arbitrary initial probability measure $ \mu \in S_ {ba} $, similar to the transformations in example 4.1. As a result (omitting the calculations) we get: $ A \mu = \frac {1} {2} \mu + \frac {1} {2} \eta $. Solve the equation $ \mu = A \mu $. From the last two equalities we obtain the only solution
  $ \mu = \eta $.

 We have shown that this combined non-degenerate MC has a unique invariant finitely additive measure $ \eta $, which is purely finitely additive, i.e. it has no nonzero countably additive component.
 
 \newpage

\section{\bf Asymptotic behavior of the norm of measures  of the decomposition components} 
 

 Consider a combined nondegenerate finitely additive MC on an arbitrary discrete space $(X, \Sigma_{d})$. 
 
 Let there be given an arbitrary initial measure $\mu^{1} \in S_{ba}$, $\mu^{1}=\mu^{1}_{ca}+\mu^{1}_{pfa}$ and  $\mu^{n+1}=A\mu^{n}, n \in N$, is the Markov sequence of measures generated by this initial measure. Let's write down its decomposition:
$$\mu^{n+1}=\mu^{n+1}_{ca}+\mu^{n+1}_{pfa}.$$
		
	
\begin{comment} 
The notation $\mu^{n+1}_{ca}$ can be interpreted in two ways: it can be a countably additive component of the measure $\mu^{n+1}$, i.e.  $(\mu^{n+1})_{ca}$, or it can be $ (n + 1) $-th iteration of measure $(\mu^{1}_{ca})$, i.e.  $(\mu^{1}_{ca})^{(n+1)}$. Generally speaking, these two interpretations do not coincide. Hereinafter, we mean that $\mu^{n+1}_{ca}=(\mu^{n+1})_{ca}$ and $\mu^{n+1}_{pfa}=(\mu^{n+1})_{pfa}$, for any $n \in N$.
\end{comment}
	
	
	
	
Since the operator $ A $ is isometric in the cone of positive measures, the norms of $\| \mu^{n+1} \| = \mu^{n+1}(X) = \| \mu^{1}(X) \| = \mu^{1}(X)  = 1$  for each  $n \in N$.

In this section, we consider the norms of the components $\| \mu^{n+1}_{ca} \|$ and $\| \mu^{n+1}_{pfa} \|$ for  $n \to \infty$, and study their composition.

Take the second iteration in the Markov sequence of measures $\mu^{2}=A\mu^{1}$. Let's make the appropriate transformations:

	$$\mu^{2} =\mu^{2}_{ca}+\mu^{2}_{pfa}=A\mu^{1} =(A_{ca}+A_{pfa})(\mu^{1}_{ca}+\mu^{1}_{pfa}) = $$ 
	$$=A_{ca}\mu^{1}_{ca}+A_{ca}\mu^{1}_{pfa}+A_{pfa}\mu^{1}_{ca}+A_{pfa}\mu^{1}_{pfa}.  \;\;\;\;\;\;\;\;\;\;  (5.1) $$
	
	
	
	
In the last four terms of the decomposition (5.1), the first is a countably additive measure, the third and fourth are purely finitely additive measures (see Theorem 3.5).

The second term $A_{ca}\mu^{1}_{pfa}$  can be a measure of any type. Consider two corresponding main cases - disjoint conditions $(H_{1})$ and $(H_{2})$. First $ (H_ {1}) $:

$$(H_{1}) \;\;\;   A_{ca}V_{pfa} \subset V_{ca},$$
that is, the operator $ A_{ca} $ transforms all purely finitely additive measures from $ V_{pfa} $ into countably additive measures. Markov chains satisfying this condition $ (H_ {1}) $ exist. It is easy to check that our Example 4.1 has this property.



\begin{theorem} 
	Let condition $ (H_{1}) $ be satisfied for a combined nondegenerate finitely additive Markov chain on an arbitrary discrete space $ (X, \Sigma_{d}) $. Then for any initial measure $ \mu^{1} \in S_{ba}$, for any $n \in N$, 
	$$ \| \mu^{n+1}_{ca} \| = \mu^{n+1}_{ca}(X) = q_{1}, $$
	$$ \| \mu^{n+1}_{pfa} \| = \mu^{n+1}_{pfa}(X) = q_{2}.$$
\end{theorem}

Proof. We carry out the proof by induction.
Let $ n = 1 $. Then by condition $ H_{1} $ the second term $A_{ca}\mu^{1}_{pfa}$ in decomposition (5.1) is a countably additive measure. Therefore, due to the uniqueness of the decomposition of the Aleksandroff-Yosida-Hewitt measures, we have:
$$\mu^{2}_{ca} = A_{ca}\mu^{1}_{ca} + A_{ca}\mu^{1}_{pfa} = A_{ca}(\mu^{1}_{ca} + \mu^{1}_{pfa}) =  A_{ca}\mu^{1}.$$
From here  $$\| \mu^{2}_{ca} \| = \mu^{2}_{ca}(X) = A_{ca}\mu^{1}(X) = q_{1} \cdot \mu^{1}(X) =q_{1}. $$
Since $$1= \| \mu^{2}\| = \mu^{2}_{ca}(X) + \mu^{2}_{pfa}(X) = \| \mu^{2}_{ca} \| + \| \mu^{2}_{pfa} \|, $$
then $$\| \mu^{2}_{pfa} \| = 1- \| \mu^{2}_{ca} \| = 1 -  q_{1} = q]_{2}.$$
Thus, the statement of the theorem for $ n = 1 $ is proved.






Suppose that the statement of the theorem is also true for some $n \in N$. 

Let's make the decomposition similar to the decomposition (5.1) for $\mu_{n+1}$ and obtain the following equalities:
$$\mu^{n+1} =\mu^{n+1}_{ca}+A\mu^{n+1}_{pfa}=A\mu^{n} =(A_{ca}+A_{pfa})(\mu^{n}_{ca}+\mu^{n}_{pfa}) =$$ $$=A_{ca}\mu^{n}_{ca}+A_{ca}\mu^{n}_{pfa}+A_{pfa}\mu^{n}_{ca}+A_{pfa}\mu^{n}_{pfa}.  \;\;\;\;\;\;\;\;\;\; \;\;\;\;\;\; \;\;\;\;\;\; (5.2) $$

As in the decomposition (5.1), here the first term is a countably additive measure, and the third and fourth terms are purely finitely additive measures.

By condition $ (H_{1}) $ the second term $A_{ca}\mu^{n}_{pfa}$ in (5.2)  is a countably additive measure. Therefore, just as for the measure $\mu^{2}_{ca}$, we obtain that $\mu^{n+1}_{ca} = A_{ca}\mu^{n}$.
In the same way, we get that
$$\| \mu^{n+1}_{ca} \| = \mu^{n+1}_{ca}(X) = A_{ca}\mu^{n}(X) = q_{1} \cdot \mu^{n}(X) =q_{1}, $$
$$\| \mu^{n+1}_{pfa} \| = \mu^{n+1}_{pfa}(X) = A_{pfa}\mu^{n}(X) = q_{2} \cdot \mu^{n}(X) =q_{2}. $$

Therefore, the statement of the theorem is true for any $ n \in N $. The theorem is proved.





\begin{comment} 
Norms $ \| \mu^{n+1}_{ca} \|$ and  $ \| \mu^{n+1}_{pfa} \|$  in Theorem 5.1 are independent of the norms of the components of the initial measure  $ \| \mu^{1}_{ca} \|$ and  $ \| \mu^{1}_{pfa} \|$. And this is a very interesting fact.
\end{comment}

\begin{corollary} 
	Let the conditions of Theorem 5.1 be satisfied. Then for such a Markov chain there exist invariant finitely additive measures $\mu^{*}=A\mu^{*}$, $\mu^{*}=\mu^{*}_{ca}+\mu^{*}_{pfa}$, and for all such measures for their components the equalities are true:
	$$ \| \mu^{*}_{ca} \| = \mu^{*}_{ca}(X) = q_{1},  \| \mu^{*}_{pfa} \| = \mu^{*}_{pfa}(X) = q_{2}.$$
\end{corollary}

Since Markov chains satisfying the condition $ (H_{1}) $ are not degenerate, that is, $0 < q_{1}, q_{2} < 1$, they do not have invariant countably additive and invariant purely finite additive measures.

Corollary 5.1 clarifies our Theorem 4.3 under the additional condition $(H_{1})$.\\



We now give the second condition $ (H_{2}) $ related to the decomposition in Theorem 5.1.
$$(H_{2}) \;\;\;   A_{ca}V_{pfa} \subset V_{pfa},$$
that is, the operator $ A_{ca} $ transforms all purely finitely additive measures from $ V_{pfa} $ into purely finitely additive measures. Such Markov chains exist, the Markov chain in Example 4.2 has this property.

\begin{theorem} 
	Let condition $ (H_{2}) $ be satisfied for a combined nondegenerate finitely additive Markov chain on an arbitrary discrete space $ (X, \Sigma_{d}) $. Then, for any initial finitely additive measure $ \mu^{1} \in S_{ba}$ for any $n \in N$ 
	$$ \| \mu^{n+1}_{ca} \| = \mu^{n+1}_{ca}(X) = q_{1}^{n} \cdot \mu^{1}_{ca}(X) = q_{1}^{n} \cdot \| \mu^{1}_{ca} \|, $$
	$$ \| \mu^{n+1}_{pfa} \| = \mu^{n+1}_{pfa}(X) = 1 - q_{1}^{n} \cdot \| \mu^{1}_{ca} \|.$$
\end{theorem}





Proof. Let us return to the decomposition (5.1). It follows from condition $ (H_{2}) $ that the second term $ A_{ca} \mu^{1}_{pfa} $ in expansion (5.1) is a purely finitely additive measure. Therefore, due to the uniqueness of the Yosida-Hewitt decomposition, $$\mu^{2}_{ca} = A_{ca}\mu^{1}_{ca}, $$
$$\mu^{2}_{pfa} = A_{ca}\mu^{1}_{pfa}+A_{pfa}\mu^{1}_{ca}+A_{pfa}\mu^{1}_{pfa} = A_{ca}\mu^{1}_{pfa}+A_{pfa}\mu^{1}. \;\;\;\;\;\;\; (5.3) $$   

Find the norm of the measure
$\mu^{2}_{ca}$ in formulas (5.3)
$$\| \mu^{2}_{ca} \| = \mu^{2}_{ca}(X) = A_{ca}\mu^{1}_{ca}(X)  = \int_{X} P_{ca}{(x,X)}\mu^{1}_{ca}(dx) =q_{1} \cdot \mu^{1}_{ca}(X) =q_{1} \cdot \| \mu^{1}_{ca} \|.$$

Since
$1= \| \mu^{2} \| =\| \mu^{2}_{ca} \| + \| \mu^{2}_{pfa} \|,$
then $\| \mu^{2}_{pfa} \|= 1 - q_{1} \cdot \| \mu^{1}_{ca} \|.$

From the equalities obtained for $ n = 1 $ $ (n + 1 = 2) $ it is still difficult to make an assumption about the general form of the norms of the components of Markov measures. Therefore, we now consider another case  $ n=2$ $(n+1=3). $





Let us make transformations for the measure $\mu^{3}$, similar to transformations (5.1) for the measure $\mu^{2}$, relying on the condition $ (H_{2}) $. As a result, we obtain equality for the measure
$$\mu^{3}_{ca} = A_{ca}\mu^{2}_{ca}$$ 
and the equality for the norm of this measure
$$\| \mu^{3}_{ca} \|= q_{1} \cdot \mu^{2}_{ca}(X)=q_{1}^{2} \cdot \mu^{1}_{ca}(X)=q_{1}^{2} \cdot \| \mu^{1}_{ca} \|.$$

From here
$$\| \mu^{3}_{pfa} \|= 1 - q_{1}^{2} \cdot \| \mu^{1}_{ca} \|.$$

Suppose now that for arbitrary $ n \in N, n \ge 2 $ holds for measures $ \mu^{n}_{ca} = A_{ca} \mu^{n-1}_{ca} $ and for the norms of these measures we have $\| \mu^{n}_{ca} \|= q_{1}^{n-1} \cdot \| \mu^{1}_{ca} \|$.

Then (omitting transformations) we get
$$\mu^{n+1}_{ca} = A_{ca}\mu^{n}_{ca},$$ 
$$\| \mu^{n+1}_{ca} \|= q_{1}^{n} \cdot \| \mu^{1}_{ca} \|,$$
$$\| \mu^{n+1}_{pfa} \|= 1 - q_{1}^{n} \cdot \| \mu^{1}_{ca} \|.$$
The theorem is proved.


\begin{comment} 
Unlike Theorem 5.1, in Theorem 5.2 the norms of the components $ \| \mu^{n+1}_{ca} \|$  and $ \| \mu^{n+1}_{pfa} \|$  of the measure   $\mu^{n+1}$ depend (linearly) on the norms of the components of the initial measure $\mu^{1}$.
\end{comment}




\begin{corollary} 
	Let the conditions of Theorem 5.2 be satisfied. Then for any finitely additive initial measure $\mu^{1} \in S_{ba}$ 
	for the components of the Markov sequence of measures generated by it $\mu^{n+1}=A\mu^{n}$ as $n \to \infty$, $\| \mu^{n}_{ca} \| \to 0 $ and   $\| \mu^{n}_{pfa} \| \to 1 $. Moreover, the convergence is uniform with respect to the initial measures $\mu^{1} \in S_{ba}$ and exponentially fast.
\end{corollary}

\begin{corollary} 
	Let the conditions of Theorem 5.2 be satisfied. Then, for such a Markov chain, all of its invariant finitely additive measures
	(and such ones always exist - see Theorem 4.1) are purely finitely additive, that is,$ \Delta_{ba} =  \Delta_{pfa} \neq \O, \Delta_{ca} = \O $. 
\end{corollary}

This statement follows from Theorem 5.2 or from its Corollary 5.2, if we take as the initial measure $\mu^{1} $ its invariant measure $\mu^{*}=A\mu^{*}$.


\end{document}